\documentclass[12pt]{article}
\usepackage{setspace}
\usepackage[margin=1in]{geometry}
\usepackage{amsmath,amssymb}
\usepackage{graphicx,multirow}
\usepackage{algorithm,algorithmic}
\usepackage{adjustbox}
\usepackage{subfigure}
\usepackage{multicol,float}
\usepackage{bbm}
\usepackage[hyphens]{url}
\usepackage{tikz,pgfplots}
\usetikzlibrary{matrix}
\usepgfplotslibrary{groupplots}
\pgfplotsset{compat=newest}
\usepackage{amsthm}
\usepackage{authblk}
\usepackage[numbers,sort]{natbib}
\usepackage[pdfborder={0 0 0}]{hyperref}
\usepackage{breakurl}
\newtheorem{theorem}{Theorem}

\newtheorem{property}[theorem]{Property}

\def\ZZ{\mathbb{Z}}

\def\RR{\mathbb{R}}

\def\N{{\cal N}}

\begin{document}

\title{Charging station optimization for balanced electric car sharing}

\author[1]{Antoine Deza\thanks{deza@mcmaster.ca}}
\author[2]{Kai Huang\thanks{khuang@mcmaster.ca}}
\author[3]{Michael R. Metel\thanks{michaelros.metel@riken.jp}}

\affil[1]{Advanced Optimization Laboratory, Department of Computing and Software, McMaster University, Hamilton, ON, Canada}
\affil[2]{DeGroote School of Business, McMaster University, Hamilton, ON, Canada}
\affil[3]{RIKEN AIP, Continuous Optimization Team, Tokyo, Japan}

\maketitle

\begin{abstract}
This work focuses on finding optimal locations for charging stations for one-way electric car 
sharing programs. The relocation of vehicles by a service staff is generally required in vehicle 
sharing programs in order to correct imbalances in the network. We seek to limit the need for 
vehicle relocation by strategically locating charging stations given estimates of traffic flow. A 
mixed-integer linear programming formulation is presented with a large number of potential 
charging station locations. A column generation approach is used which finds an optimal set of 
locations for the continuous relaxation of our problem. Results of a numerical experiment using 
real traffic and geographic information system location data show that our formulation significantly 
increases the balanced flow across the network, while our column generation technique 
was found to produce a superior solution in much shorter computation time compared to solving the 
original formulation with all possible station locations.
\end{abstract}

{\bf Keywords:} electric vehicle; one-way car sharing; column generation; mixed-integer linear 
programming

\section{Introduction}

Electric car sharing programs are a method for urban centres to combat traffic congestion and 
pollution \citep{boyaci2017,brand2017} as well as to promote the use of green technologies. In 
one-way electric car sharing programs such as \citet{autolib} in Paris, France, users are able to 
use and return vehicles to any charging station in the network. Large imbalances with the supply of 
vehicles and parking spaces across nodes in the network are generally observed, requiring a service 
staff to continuously transport vehicles to satisfy demand.\\ 

The problem of determining the optimal locations of charging stations for 
electric car sharing systems is considered in \citep{brand2017}. The number of stations and 
vehicles, and their optimal placement is determined in order to maximize profit. One assumption is 
that decision makers do not consider operational activities of a service staff, and in particular 
for vehicle relocation, though cars must be initially placed at stations at the beginning of the 
optimization time horizon. In the scenario where an electric car sharing network has already been 
built, a number of researchers have developed methodologies for the vehicle relocation problem. The 
use of folding bicycles by workers which fit into the trunks of cars \citep{brug2014}, and the use 
of towtrucks which are capable of moving a number of vehicles at a time \citep{dror1998} have been 
proposed. In \cite{ait2018}, a set of agents are assumed to be employed which drive vehicles 
between 
stations. Given a set of predetermined trips, the number of vehicles, agents, and the schedule of 
the agents are optimized each day. In addition, the optimal relocation of workers themselves across 
the network in order to relocate vehicles has been considered in \cite{boyaci2017}.\\

The idea of having to relocate vehicles runs counter to the objective of decreasing traffic 
congestion, and will cut into profit and system efficiency. The inefficiency of transporting 
vehicles between stations is compounded for electric vehicles, as not only are they not being 
used productively while being transported, but will require further charging afterwards if driven 
by workers between stations.\\

In this work we consider a set of nodes with an expected traffic flow between each pair of nodes per 
time period. A methodology is presented to place charging stations in such 
a way so as to limit supply imbalances in the network by matching demand for vehicles and parking 
spaces at each charging station. By a strategic placement of charging stations, the need for a 
service staff to continuously relocate vehicles is greatly reduced, which is shown by a 
significant increase in the estimated balanced flow across all charging stations.\\

The remainder of the paper is organized as follows. In Section 2 we describe the problem setting, 
the concept of a balanced electric car sharing network, and its corresponding optimization 
problem (\ref{opt1}). Section 3 is where the balanced charging station algorithm 
(BCSA) is developed using a column generation technique to find a near optimal set of charging 
stations. Section 4 outlines how to solve (\ref{opt1}) by exhaustively enumerating all possible 
charging stations which is used for comparison with BCSA. In Section 5 we describe the numerical 
experiment with the use of real traffic and geographic information system location data, 
and present the results. The paper finishes in Section 6 with the conclusion.\\  

\section{Balanced electric car sharing optimization model}

Let $N$ be the set of trip nodes, which is the set of locations where trips are arriving to 
and departing from, and let $T=\{1,2,...,M_T\}$ be a set of time intervals over a 24 hour cycle 
with lengths $L_t$. For each $t\in T$ there exists an origin-destination matrix $OD^t\in 
\ZZ^{|N|\times |N|}$, indicating for each pair of nodes $\{n,n'\}\in N$, an estimate of the 
number of trips from $n$ to $n'$. For each node $n\in N$, its outward flow over $t$ is 
$f^{t-}_n=\sum_{n'\in N}OD^t_{n,n'}$, which is the sum of trips departing from node $n$, 
requiring an electric vehicle near $n$. Its inward flow over $t$ is $f^{t+}_n=\sum_{n'\in 
N}OD^t_{n',n}$, which is the sum of trips arriving to node $n$, requiring a parking 
space near $n$.\\

Let $S$ be the set of potential charging stations. We assume that people are willing to walk up to 
$w=0.5$ km to or from a charging station as used in \citep{ait2018, brug2014}. For each $s\in S$ we 
define its set of neighbouring nodes as a subset of nodes which are within walking distance from it, 
$\N(s)\subseteq\{n\in N: d(n,s)\leq w\}$, where $d(n,s)$ is the distance in kilometres between $n$ and 
$s$. Likewise, we define the neighbouring set of stations for each $n \in N$ as a 
subset $\N(n)\subseteq\{s\in S: d(s,n)\leq w\}$. In addition, $n\in \N(s)$ if and only if $s\in 
\N(n)$. We assume all distances are calculated using the $l_1$-norm.\\

We must decide how many parking spaces to place at each charging station. As will be described in 
more detail, we require a balanced flow at each station, meaning the number of trip arrivals and 
departures assigned to a station are always equal. Given this requirement, assigned trips will be 
arriving and departing with probability 0.5, so we always choose an even number of parking 
spaces at each charging station, with the intention of having half the number of vehicles as spaces 
to best meet incoming demand. This is in line with what is observed in current use, with 2 and 4 
parking spaces being the most common for a charging station \citep{frade2011}.\\

In order to determine the usage capacity of a pair of parking spaces, we need an estimate of 
how long it will take for someone to park and plug in a vehicle, and to register a vehicle 
and leave the station. Let $p_a$ and $p_d$ be the times to perform these two tasks in hours, 
respectively. As will be clear, these individual estimates are not required, but only their 
average $p=\frac{1}{2}(p_a+p_d)$. We also need an estimate of the average trip length, from which 
we can estimate the average charging time required after each trip. For each time period $t$, we 
calculate $l_t$, the average trip length over the network, 

$$l_{t}=\frac{\sum_{n\in N}\sum_{n'\in N}OD^t_{n,n'}d(n,n')}
{\sum_{n\in N}\sum_{n'\in N}OD^t_{n,n'}}.$$ 

Given an estimate of the electric vehicle's charging time per kilometre driven, $u$ in hours/km, we 
can estimate the average 
charging time required after each trip, $a_t=ul_t$. The amount of time on average 
required between trips for a car to be dropped off and recharged, or picked up is then
$p+\frac{a_t}{2}$. We note that this average time is valid as the number of assigned arrivals and 
departures at each station are always equal. Assuming that $OD$ contains an estimate of the total 
traffic flow within a city, we assume that $k=0.5\%$ of the flow will eventually be fulfilled using 
our service as in \citep{brug2014}. The maximum amount of flow a pair of parking spaces at station $s$ 
can be allocated during time period $t\in T$ is then estimated as  

$$v_{t}=2\left\lfloor\frac{L_t}{k(p+\frac{a_t}{2})}\right\rfloor$$

where the estimated capacity of a single parking space, $L_t/(k(p+a_t/2))$, has been rounded down to 
the nearest integer to determine the maximum number of whole trips. We note that 
$v_t$ is an even integer, which results in integral solutions of our optimization problem 
(\ref{opt1}), as proved in Property \ref{prop1}.\\
 
Whereas the objective of a company running a combustion engine car sharing program is most likely 
to be to maximize profit, in this work we assume that the electric car sharing 
program is funded or heavily subsidized by a government in order to reduce pollution and road 
congestion, and to promote green technologies. For this reason, the objective in our 
optimization problem (\ref{opt1}) is to maximize the number of electric vehicle trips, written 
equivalently as minimizing the number of trips not satisfied. In the literature, our objective is 
most similar to that of \cite{ait2018} and \cite{boyaci2017}, which both used multi-objective 
models, minimizing costs and unsatisfied customers. We as well minimize the number of unsatisfied 
customers, but place costs in a budget constraint.\\

The flow to and from each node is assigned to its neighbouring stations. $F^{t-}\in \RR^{|S|\times 
|N|}_+$ is the decision matrix of the number of trips from each node $n$ assigned to leave from 
each station $s$ during time period $t$, and likewise $F^{t+}\in \RR^{|S|\times |N|}_+$ is the 
decision matrix of the number of trips to each node $n$ assigned to arrive at each station $s$ 
during time period $t$. Further, $E^t\in \RR^{|N|\times |N|}_+$ is the traffic flow in $OD^t$ 
which is not satisfied by any station, and $e$ is a vector of ones of size $|N|$. We assign the 
inward flow of nodes to their neighbouring stations by constraint set (2), and their outward flow 
by (3). Care must be taken to maintain the integrity of each trip, meaning if an incoming trip to 
station $n$ is not assigned to a charging station, then a corresponding trip leaving from a node 
$n'$ to $n$ cannot be assigned to a station's outward flow. If an inward trip to node $n$ is 
not satisfied in time $t$, then an entry $n'$ of the $n^{th}$ column of $E^t$ is incremented by 1 in 
(2). This then forces one corresponding outward trip from $n'$ to not be satisfied in (3). This 
is further enforced by constraint (5), which ensures that the entry in $E^t$ corresponds to an 
actual trip in $OD$. A small example showing the necessity of constraint (5) is found in the 
Appendix subsection titled ``Requirement of (5) in (\ref{opt1}) for trip integrity".\\   

$z_s$ is the decision variable for the number of pairs of parking spaces to 
install at station $s$, $m_s$ is the maximum number of pairs of parking spaces that can be 
installed, and $c_s$ is the cost of constructing a pair of parking spaces. We assume all costs, 
including the electric vehicles themselves are embedded into the price per pair of parking 
spaces. $b$ is the budget allocated for the construction of the electric car sharing network. 
Constraint set (4) enforces the net flow over each station to be zero, so as to minimize network 
imbalances and minimize the need for vehicle relocation, as the expected number of cars parking 
and leaving at each charging station over each time period are equal. We have also included the 
corresponding dual 
variables to the left of each set of constraints. 

\begin{alignat}{6}\label{opt1}\tag{OP}
&&\min&\text{ }\sum_{t\in T}\sum_{n,n'\in N}E^t_{n,n'}&\\
\left(U_{t,s}\geq 0\right)\hspace{10 pt}&&\mbox{s.t. }&\sum_{n\in 
\N(s)}\left(F^{t+}_{s,n}+F^{t-}_{s,n}\right)\leq 
v_tz_s\hspace{10 pt}\forall\text{ }t\in T\text{, }s\in S&&\left(1\right)\nonumber\\
\left(P_{t,n}\right)\hspace{10 pt}&&&\sum_{s\in \N(n)} F^{t+}_{s,n} + 
(E^t_{\cdot,n})^Te=f^{t+}_n\hspace{10 
pt}\forall\text{ }t\in T\text{, }n\in N&&\left(2\right)\nonumber\\
\left(G_{t,n}\right)\hspace{10 pt}&&&\sum_{s\in \N(n)} F^{t-}_{s,n} + 
E^t_{n,\cdot}e=f^{t-}_n\hspace{10 pt}\forall\text{ }t\in T\text{, }n\in N&&\left(3\right)\nonumber\\
\left(R_{t,s}\right)\hspace{10 pt}&&&\sum_{n\in 
\N(s)}\left(F^{t+}_{s,n}-F^{t-}_{s,n}\right)=0\hspace{10 
pt}\forall\text{ }t\in T\text{, }s\in S&&\left(4\right)\nonumber\\
\left(W^t_{n,n}\geq 0\right)\hspace{10 pt}&&&E\leq OD&&\left(5\right)\nonumber\\
\left(q\geq 0\right)\hspace{10 pt}&&&\sum_{s\in S} c_sz_s\leq b&&\nonumber\\
\left(h_s\geq 
0\right)\hspace{10 pt}&&&z_s\leq m_s\hspace{10 pt}\forall\text{ }s\in S&&\nonumber\\
&&&F^{t+}\in \RR^{|S|\times |N|}_+,\hspace{10 pt} F^{t-}\in \RR^{|S|\times 
|N|}_+,\hspace{10 pt} E^t\in \RR^{|N|\times |N|}_+, \hspace{10 pt} z\in \ZZ^{|S|}_+\nonumber
\end{alignat}

\begin{property}\label{prop1}	
	The optimal solution of (\ref{opt1}) is integral.	
\end{property}

\begin{proof}

Given that $F^{t+}_{s,n}=F^{t-}_{s,n}$ from (4), the station capacity constraints (1) can be written 
equivalently as 

\begin{alignat}{6}
\hspace{95 pt}&&&\sum_{n\in \N(s)}F^{t+}_{s,n}\leq \frac{v_t}{2}z_s\hspace{10 pt}\forall\text{ 
}t\in 
T\text{, }s\in S&&\hspace{125 pt}\left(1'\right)\nonumber
\end{alignat}

where the right-hand side, $v_tz_s/2\in \ZZ$. Let $F_{\text{vec}}^{t+}$ be the vector of all columns 
of $F^{t+}$ stacked on top of 
each other, where all elements corresponding to an $s\notin\N(n) \text{ }\forall n$ have been 
removed. 
Let $F_{\text{vec}}^{t-}$ be built in the same manner, and let $E^t_{\text{vec}}$ simply be the 
vector of all columns of $E^{t}$ stacked on top of each other. We set 
$X^t=[F^{t+}_{vec};F^{t-}_{vec};E^t_{vec}]$ and can now write constraints [(1'),(2)-(5)] as 
$AX^t\{\leq,=\} b^t$, where $A$ is a matrix composed of $\{-1,0,1\}$, $b^t$ is a vector composed 
of integers, and $\{\leq,=\}$ stands in for $\leq$ or $=$ as appropriate for each row. As a 
reference, we have included a subsection of the Appendix entitled ``An example of $AX^t\{\leq,=\} 
b^t$". To show that A is totally unimodular we use Property 2 which is borrowed from 
\citep{schrijver1998}.

\begin{property}{\citep[Theorem 19.3]{schrijver1998}}
	Let $A$ be a matrix with entries 0, 1, or -1. $A$ is totally unimodular if and only if each 
	collection $J$ of columns of $A$ can be split into two parts, $J_1$ and $J_2$, so that 
	the sum of the columns in $J_1$ minus the sum of the columns $J_2$ is a vector with entries 
	only equal to 0, 1, or -1.			
\end{property}

As the transpose of a totally unimodular matrix is also totally unimodular, we show that the 
condition holds for each collection of rows of $A$. We will assign rows from each constraint set to 
$J_1$ or $J_2$, and show that for each column the difference of the sums has magnitude 1 or 0.\\

We partition the columns of $A$ into three subsets $\{C_1,C_2,C_3\}$ by which variables they 
multiply. Let $C_1$ be the first $\sum_{n\in N}|\N(n)|$ columns multiplying $F^{t+}_{vec}$, let $C_2$ 
be the second $\sum_{n\in N}|\N(n)|$ columns multiplying $F^{t-}_{vec}$, and let $C_3$ be the last  
$|N|\times|N|$ columns multiplying $E^{t}_{vec}$.\\

We begin with $C_2$, and note that each variable $F^{t-}_{s,n}$ will be found in one constraint in 
(3), with a coefficient of 1 in A, and in one constraint of (4), with a coefficient of -1. We 
therefore place all rows from (3) and (4) into $J_2$. The sum of each column over $C_2$ will then be 
either 0, 1, or -1. The current partitioning is $J_1=\{\emptyset\}$ and $J_2=\{(3),(4)\}$.\\

For $C_3$ each variable $E_{n,n'}$ is found in one constraint in (2), (3), and (5), 
all with a coefficient of 1 in $A$. As rows from (3) are in $J_2$, we place all rows from (2) in 
$J_1$. We note that rows from (5) always contain only one non-zero element equal to 1: They are zero 
over $C_1$ and $C_2$, and form an identity matrix in $C_3$. This implies that their placement in 
$J_1$ or $J_2$ only affects the sums of columns in $C_3$. Any row $i$ in (5) with 1 in column $j$ can 
then be placed in either $J_1$ or $J_2$ to ensure that column $j$ sums to 0, -1, or 1: If the rows 
from (2) and (3) with 1 in column $j$ are in $J$, or neither are in $J$, then $i$ can be 
placed in $J_1$ or $J_2$. If only the row from (2) is in $J$, then place $i$ in $J_2$, and if only 
the row from (3) is in $J$, place $i$ in $J_1$. The current partitioning is $J_1=\{(2),(5)_1\}$ and 
$J_2=\{(3),(4),(5)_2\}$, where (5) is partitioned as described.\\

For $C_1$ each variable $F^{t+}_{s,n}$ is found in one constraint in (1'), (2), and (4), all 
with a coefficient of 1. We have already placed rows from (2) in $J_1$ and rows from (4) in $J_2$. 
Given that constraints (1') and (4) both contain the same summation of $F^{t+}$ over all $t$ and $s$, 
the rows of $A$ corresponding to (1') and (4) are identical over $C_1$, meaning that for each row 
of (1') there is a copy in the rows of (4), in $C_1$. In addition, the rows of $C_1$ 
only multiply variables from $F^{t+}_{vec}$, meaning their placement in either $J_1$ or $J_2$ only 
affect the sums of columns of $C_1$. We can therefore place rows from (1') in $J_1$ or $J_2$ so 
that each column of $C_1$ sums to 0, 1, or -1 as done previously with rows of (5) for $C_2$. The 
final partitioning is $J_1=\{(1')_1,(2),(5)_1\}$ and $J_2=\{(1')_2,(3),(4),(5)_2\}$.

\end{proof}

Property \ref{prop1} would not hold in general without the fact that $v_t$ is even for all $t$, as can 
be observed in the following scenario. For an optimal solution where $v_t$ is odd for a constraint in 
(1) which is binding, half of an inflow of a trip and half of an outflow of a trip would have to be 
assigned to $s$, as $\sum_{n\in \N(s)}F^{t+}_{s,n}$ and $\sum_{n\in \N(s)}F^{t-}_{s,n}$ are equal, but 
sum to an odd number.\\

An undetermined aspect of (\ref{opt1}) is the set of stations $S$. If we were to 
consider every location in the city it would be difficult to bound its size. We limit our set of 
potential stations to one per feasible subset of nodes, so that $|N|\leq |S| \leq 2^{|N|}-1$, 
where the lower bound occurs when $d(n,n')>2w$ and the upper bound occurs when $d(n,n')\leq 2w$ 
for all $n,n'\in N$. Though we don't expect $|S|$ to be near its upper bound, its size has the 
potential to 
make (\ref{opt1}) a computationally challenging problem, particularly when implemented in dense urban 
areas with many nodes. For this reason, we propose a column generation approach outlined in the next 
section.

\section{Determining $S$ using column generation}

We begin with some initial set of stations $S'$ and 
iteratively add stations using a column generation technique until we have found an optimal 
set of stations for the continuous relaxation of (\ref{opt1}), and use this set to find a solution 
to (\ref{opt1}). The dual program (\ref{opt2}) of the continuous relaxation of (\ref{opt1}) will 
be used in the solution technique.

\begin{alignat}{6}\label{opt2}\tag{D}
&\max&&\text{ }-bq-\sum_{s\in S}m_sh_s-\sum_{t\in T}\sum_{n,n'\in N}OD^t_{n,n'}W^t_{n,n'}-\sum_{t\in 
	T}\sum_{n\in N}(f^{t+}_nP_{t,n}+f^{t-}_nG_{t,n})\nonumber\\
&\mbox{s.t. }&&1+G_{t,n}+P_{t,n'}+W^t_{n,n'}\geq 0\hspace{10 pt}\forall\text{ }t\in T\text{, }n,n'\in 
N&&\left(1\right)\nonumber\\
&&&U_{t,s}+P_{t,n}+R_{t,s}\geq 0\hspace{10 pt}\forall\text{ }t\in T\text{, }s\in S\text{, }n\in 
\N(s)&&\left(2\right)\nonumber\\
&&&U_{t,s}+G_{t,n}-R_{t,s}\geq 0\hspace{10 pt}\forall\text{ }t\in T\text{, }s\in S\text{, }n\in 
\N(s)&&\left(3\right)\nonumber\\
&&&h_s+qc_s-\sum_{t\in T}v_tU_{t,s}\geq 0\hspace{10 pt}\forall\text{ }s\in 
S&&\left(4\right)\nonumber\\
&&&U\in \RR^{M_T\times |S|}_+\nonumber, \hspace{10 pt} W\in \RR^{T\times|N|\times |N|}_+, 
\hspace{10 pt} h\in \RR^{|S|}_+,q\geq 
0\nonumber
\end{alignat}

\subsection{Finding a new station}
\label{FNS}

A new feasible station $s'$ to add to $S'$ will satisfy (\ref{opt3}). The 
neighbourhoods of all current stations are encoded in a matrix $B\in \{0,1\}^{|N|\times |S'|}$, 
where $B_{\cdot,s}$ 
represents $\N(s)$. A new column will be added for $s'$, with constraint set (2) ensuring its 
uniqueness. We define $d_n$ as the maximum feasible distance of a station from $n$, which equals 
$d_n=\max_{n'\in N} d(n,n')$. This is used as a big-M constant when ensuring each $n\in \N(s')$ is 
within $w$ of $s'$ in constraint set (1). 
Each node $n$ has coordinates $x_n$ and $y_n$, stored in vectors $x\in \RR^{|N|}$ and $y\in 
\RR^{|N|}$. Let $s'_x$ and $s'_y$ be the coordinates of $s'$, which we write as a convex 
combination of the coordinates of its neighbouring nodes using constraints (3) and (4). This 
enables us to infer the cost of installing a pair of parking spaces at $s'$, and its parking 
space capacity. Each node $n$ has an estimated cost $c^N_n$ for constructing a pair of parking 
spaces at $n$, as well as a capacity of pairs of parking spaces $m^N_n$. We 
estimate the price $c_{s'}$ and the capacity $m_{s'}$ of $s'$ using the weighted 
average of the nodes' values within $\N(s')$ in constraints (5) and (6). Constraint (6) sets 
$m_{s'}$ to the rounded value of $\alpha^Tm^N$.

\begin{alignat}{6}\label{opt3}\tag{FS}
&d(s',n)\leq w + (1-B_{n,s'})(d_n-w) \hspace{10 pt}\text{for }n\in N&&\left(1\right)\\
&B_{\cdot,s}^TB_{\cdot,s'} + (1-B_{\cdot,s})^T(1-B_{\cdot,s'})\leq |N|-1\hspace{10 pt}\text{for 
}s\in S'\nonumber&&\left(2\right)\\
&s'_x=\alpha^Tx, \hspace{10 pt}s'_y=\alpha^Ty\nonumber&&\left(3\right)\\
&\alpha^Te=1, \hspace{10 pt}\alpha\leq B_{\cdot,s'}\nonumber&&\left(4\right)\\
&c_{s'}=\alpha^Tc^N\nonumber&&\left(5\right)\\
&\alpha^Tm^N-0.5\leq m_{s'}\leq\alpha^Tm^N+0.5\nonumber&&\left(6\right)\\
&\alpha\in\RR^{|N|}_+, \hspace{10 pt} B_{n,s'}\in\{0,1\}^{|N|}, \hspace{10 pt} c_{s'}\in 
\RR_{++}, \hspace{10 pt} m_{s'}\in \ZZ_{++}\hspace{100 pt}&&\nonumber
\end{alignat}

After solving the dual program (\ref{opt2}) for some set $S'$, we seek to find a feasible station 
satisfying (\ref{opt3}) which could result in a decrease in the objective of (\ref{opt2}). The 
dual objective will decrease if $h_{s'}>0$, or equivalently if $\sum_{t\in 
T}v_tU_{t,s'}-qc_{s'}>0$ from constraint (4). From constraint sets (2) and (3) of (\ref{opt2}), it 
follows that $U_{t,s'}\geq -\frac{(P_{t,n}+G_{t,n'})}{2}$ for all $n,n'\in \N(s')$
and this inequality will be binding for some $n,n'\in \N(s')$ or equal to 0. We build the 
matrices $PG^t\in\RR^{|N|\times |N|}_+$, where 
$PG^t_{n,n'}=\max\left(-\frac{(P_{t,n}+G_{t,n'})}{2},0\right)$. The optimal value of $U_{t,s'}$ in 
(\ref{opt2}) will be the maximum value of $PG^t$ over $n,n'\in \N(s')$, holding current dual 
variables constant. We can determine $U_{t,s'}$ with (\ref{opt5}), where $PG^t\circ D^t$ is taking 
the Hadamard product of the two matrices.
\begin{alignat}{6}\label{opt5}\tag{DU}
&\max&&\text{ }\sum_{t\in T}U_{t,s'}\\
&\mbox{s.t. }&&U_{t,s'}=\sum_{n\in N}\sum_{n'\in N}(PG^t\circ D^t)_{n,n'}\hspace{10 
pt}\text{for }t\in 
T\nonumber\\
&&&e^TD^te=1\hspace{10 pt}\text{for }t\in T\nonumber\\
&&&D^t_{n,n'}\leq B_{n,s'}\hspace{10 pt}\text{for }t\in 
T,n,n' \in N\nonumber\\
&&&D^t_{n,n'}\leq B_{n',s'}\hspace{10 pt}\text{for }t\in 
T,n,n'\in N\nonumber\\
&&&D\in \RR^{|N|\times |N| \times M_T}_+\nonumber
\end{alignat}
Putting (\ref{opt3}) and (\ref{opt5}) together, we find the next charging station $s'$ to 
include in (\ref{opt1}) by solving (\ref{opt6}). If $h_{s'}>0$ we have found a station vector 
$B_{\cdot,s'}$ which can improve the optimal solution of (\ref{opt1}). 
\begin{alignat}{6}\label{opt6}\tag{NS}
&\max&&\text{ }h_{s'}=\sum_{t\in T}v_{t}U_{t,s'}-qc_{s'}\\
&\mbox{s.t. }&&U_{t,s'}=\sum_{n\in N}\sum_{n'\in N}(PG^t\circ D^t)_{n,n'}\hspace{10	pt}\text{for 
}t\in T\nonumber\\
&&&e^TD^te=1\hspace{10 pt}\text{for }t\in T\nonumber\\
&&&D^t_{n,n'}\leq B_{n,s'}\hspace{10 pt}\text{for }t\in 
T,n,n'\in N\nonumber\\
&&&D^t_{n,n'}\leq B_{n',s'}\hspace{10 pt}\text{for }t\in 
T,n,n'\in N\nonumber\\
&&&d(s',n)\leq w + (1-B_{n,s'})(d_n-w) \hspace{10 pt}\text{for }n\in N\nonumber\\
&&&B_{\cdot,s}^TB_{\cdot,s'} + (1-B_{\cdot,s})^T(1-B_{\cdot,s'})\leq N-1\hspace{10 pt}\text{for 
}s\in S'\nonumber\\
&&&s'_x=\alpha^Tx,\hspace{10 pt}s'_y=\alpha^Ty\nonumber\\
&&&\alpha^Te=1, \hspace{10 pt}\alpha\leq B_{\cdot,s'}\nonumber\\
&&&c_{s'}=\alpha^Tc^N,\hspace{10 pt}\alpha\in\RR^{|N|}_+\nonumber\\
&&&B_{\cdot,s'}\in\{0,1\}^{|N|},\hspace{10 pt}D\in 
\RR^{|N|\times |N|\times M_T}_+\nonumber
\end{alignat}

We did not consider the effect of $m_{s'}$ in (\ref{opt6}), and so we must refine the station's 
location. Let $\overline{f}^t_{s'}:=2\min(B^T_{\cdot,s'}f^{t-},B^T_{\cdot,s'}f^{t+})$, 
which is the maximum flow that could be allocated to station $s'$ in time $t$ if all of its 
neighbouring nodes allocate all of their trips to $s'$. Let $m^c_{s'}:=\max\limits_{t\in 	
T}\left\lceil\frac{\overline{f}^t_{s'}}{v_t}\right\rceil$, which is the maximum parking space 
capacity station $s'$ could require. 
In (\ref{opt7}), we optimize $\alpha$ given our optimal solution for  $B_{\cdot,s'}$ and 
$U_{t,s'}$ from (\ref{opt6}). We do not 
want to reward a location for capacity which will not be used, so we use $\min(m_{s'},m^c_{s'})$ in 
the objective. 

\begin{alignat}{6}\label{opt7}\tag{OA}
&\max&&\text{ }\min(m_{s'},m^c_{s'})\left(\sum_{t\in T}v_{t}U_{t,s'}-qc_{s'}\right)\\
&\mbox{s.t. }&&s'_x=\alpha^Tx,\hspace{10 pt}s'_y=\alpha^Ty\nonumber\\
&&&\alpha^Te=1, \hspace{10 pt}\alpha\leq B_{\cdot,s'}\nonumber\\
&&&c_{s'}=\alpha^Tc^N,\hspace{10 pt}m_{s'}\leq \alpha^Tm^N+0.5\nonumber\\
&&&\alpha\in\RR^{|N|}_+,\hspace{10 pt}m_{s'}\in\ZZ_{++}\nonumber
\end{alignat}

As a non-convex problem, we solve (\ref{opt7}) as a set of linear programs, (\ref{opt8}), 
for $m^i_{s'}$ over the 
interval   
$$m^{\text{min}}_{s'}\leq m^i_{s'}\leq 
\max\left(m^{\text{min}}_{s'},m^{\bar{c}}_{s'}\right).$$ 
where $m^{\text{min}}_{s'}:=\min\limits_{n\in \N(s')}m_n$, 
$m^{\bar{c}}_{s'}:=\min\left(m^c_{s'},m^{\text{max}}_{s'}\right)$, and 
$m^{\text{max}}_{s'}:=\max\limits_{n\in \N(s')}m_n$. We take the optimal $\alpha$ vector from the 
program (\ref{opt8}) with the maximum objective.

\begin{alignat}{6}\label{opt8}\tag{$\text{OA}^i$}
&\max&&\text{ }m^i_{s'}\left(\sum_{t\in T}v_{t}U_{t,s'}-qc_{s'}\right)\\
&\mbox{s.t. }&&s'_x=\alpha^Tx,\hspace{10 pt}s'_y=\alpha^Ty\nonumber\\
&&&\alpha^Te=1, \hspace{10 pt}\alpha\leq B_{\cdot,s'}\nonumber\\
&&&c_{s'}=\alpha^Tc^N,\hspace{10 pt}m^i_{s'}\leq \alpha^Tm^N+0.5\nonumber\\
&&&\alpha\in\RR^{|N|}_+\nonumber
\end{alignat}

\subsection{Valid inequalities for (\ref{opt6})}

Many instances of (\ref{opt6}) must be solved, which can be time consuming given the binary 
variables $B_{\cdot,s'}$. A set of valid inequalities were added to (\ref{opt6}) which were found 
to significantly reduce computation time. The general idea is that if we 
are 
given a subset of nodes $N'$ where $d(n,n')>2w$ for all $n,n'\in N'$, then 
$\sum_{n\in N'}B_{n,s'}\leq 1$. The algorithm for adding these constraints is found in Algorithm 
\ref{alg:alg1}. 

\begin{algorithm}[H]
	\caption{Valid inequalities for (\ref{opt6})}
	\label{alg:alg1}
	\begin{algorithmic}[1]
		\FOR{$n\in N$}
		\STATE{$V=\{n\}$}
		\FOR{$n'\in N$}
		\IF{$d(y,n')>2w$ $\forall y\in V$}
		\STATE{$V=V+\{n'\}$}
		\ENDIF
		\ENDFOR
		\IF{$|V|>1$}
		\STATE{add $\sum_{y\in V}B_{y,s'}\leq 1$ to  (\ref{opt6})}
		\ENDIF
		\ENDFOR		
	\end{algorithmic}
\end{algorithm}

\subsection{Column generation search heuristic}

We write the location of $s'$ as a convex combination of nodes as a means of estimating its cost 
and capacity, but to motivate the proposed heuristic, let us first assume $s_x'$ and $s_y'$ are 
free variables, and the station's cost and capacity are not functions of its neighbouring nodes.\\

If we consider the problem setting where $M_T=1$ then an optimal solution for the convex relaxation 
of (\ref{opt1}) can be found consisting of stations of size no greater than 2. This is due to 
when maximizing the objective of (\ref{opt6}), $PG$ is only a function of two nodes. More 
intuitively, given a number of nodes all within $2w$ of eachother, it will always be optimal to 
build a number of fractional stations for each pair of nodes, which gives maximum flexibility of 
location. 
Increasing the size of $M_T$, a neighbourhood of size $|\N(s')|\leq 2M_T$ will be optimal, where 
$|\N(s')|=2M_T$ can occur when for each $t$, new nodes make up the pair mapping to the maximum of 
$PG^t$.\\

In our setting, this bound no longer applies as adding extra nodes can decrease station cost and 
increase capacity, but when our technique creates stations with larger neighbourhoods it is not for 
our objective of matching flow from different nodes, but for what can be considered secondary 
concerns relating to our constraints. In order to counter this problem, in each iteration we add a 
station with the largest possible neighbourhood which still results in a positive objective value 
for (\ref{opt6}). To begin, we solve (\ref{opt9}) to 
find the largest possible neighbourhood size. We then proceed to solve (\ref{opt6}) with the added 
constraint $\sum_{n\in N}B_{n,s'}=\overline{SB}$ which we will refer to as (NS($\overline{SB}$)). 
The value of $\overline{SB}$ is then decremented when no more stations can be found which result in 
$h_{s'}>0$. Details of this process are found in the Subsection \ref{BCSA}.

 \begin{alignat}{6}\label{opt9}\tag{SB}
 &\max&&\text{ }\overline{SB}:=\sum_{n\in N} B_{n,s'}\\
 &\mbox{s.t. }&&d(s',n)\leq w + (1-B_{n,s'})(d_n-w) \hspace{10 pt}\text{for }n\in N\nonumber\\
 &&&s'_x=\alpha^Tx, \hspace{10 pt}s'_y=\alpha^Ty\nonumber\\
 &&&\alpha^Te=1, \hspace{10 pt}\alpha\leq B_{\cdot,s'}\nonumber\\
 &&&\alpha\in\RR^{|N|}_+,\hspace{10 pt} B_{n,s'}\in\{0,1\}^{|N|}\nonumber
 \end{alignat}
 
 \subsection{$S'$ initialization}
 
 We initialize $S'$ as the station which would satisfy the greatest balanced flow, which is found 
 by solving (\ref{opt11}).  
 
 \begin{alignat}{6}\label{opt11}\tag{GF}
 &\max&&\text{ }\sum_{t\in T}\min(2B^Tf^{t-},2B^Tf^{t+},m_{s'}v_t)\\
 &\mbox{s.t. }&&d(s',n)\leq w + (1-B)(d_n-w) \hspace{10 pt}\text{for }n\in N\nonumber\\
 &&&s'_x=\alpha^Tx,\hspace{10 pt}s'_y=\alpha^Ty\nonumber\\
 &&&\alpha^Te=1, \hspace{10 pt}\alpha\leq B\nonumber\\
 &&&m_{s'}\leq \alpha^Tm^N+0.5,\hspace{10 pt}m_{s'}\in\ZZ_{++}\nonumber\\
 &&&\alpha\in\RR^{|N|}_+,\hspace{10 pt} B\in\{0,1\}^{|N|}\nonumber
 \end{alignat}
 
If $m^{\text{max}}_{s'}-m^{\text{min}}_{s'}>0$  there could be a trade-off between capacity and 
price when determining the station's location. We minimize the station price without restricting 
potential flow by solving (\ref{opt12}), where $m^{\bar{c}}_{s'}$ is as defined in Subsection 
\ref{FNS}, and we set $m_{s'}=\lceil\alpha^Tm^N\rfloor$. We call the station found from this 
process 
$s^{GF}$.

 \begin{alignat}{6}\label{opt12}\tag{$OA^{\bar{c}}$}
 &\min&&\text{ }\alpha^Tc^N\\
 &\mbox{s.t. }&&s'_x=\alpha^Tx,\hspace{10 pt}s'_y=\alpha^Ty\nonumber\\
 &&&\alpha^Te=1, \hspace{10 pt}\alpha\leq B_{\cdot,s'}\nonumber\\
 &&&m^{\bar{c}}_{s'}\leq \alpha^Tm^N+0.5\nonumber\\
 &&&\alpha\in\RR^{|N|}_+\nonumber
 \end{alignat}

\subsection{Balanced car sharing algorithm}
\label{BCSA}
We now present the BCSA algorithm for finding a balanced electric car sharing charging station 
network. 

\begin{algorithm}[H]
	\caption{Balanced charging station algorithm (BCSA)}
	\label{alg:alg2}
	\begin{algorithmic}[1]
		\STATE{Initialize $S'=s^{GF}$}
		\STATE{Solve (\ref{opt9})}
		\STATE{Solve (\ref{opt2})}\label{dual}
		\WHILE{$\overline{SB}>0$}
		\STATE{Solve (\ref{opt6}($\overline{SB}$))}
		\IF{$h_{s'}>0$}
		\STATE{Solve (\ref{opt7})}
		\STATE{$S'=S'\cup \{s'\}$} 
		\STATE{Solve (\ref{opt2})}\label{dual}
		\ELSE
		\STATE{$\overline{SB}=\overline{SB}-1$}
		\ENDIF		
		\ENDWHILE
		\STATE{Solve (\ref{opt1}) with $S=S'$}
	\end{algorithmic}
\end{algorithm}

\section{Exhaustive enumeration method}

Instead of iteratively finding stations, we consider finding all possible stations $S$ initially 
and then solving (\ref{opt1}) directly. We again limited the search to a single station per subset 
of 
nodes. The basic means of finding $S$ was as follows. For each node $n$, we found the set of nodes 
$S_n=\{n':d(n,n')\leq 2d,n'\neq n\}$. We then found all subsets $S'_n$ of size  
$|S_n|$,$|S_n|-1$,...,$0$ in $S_n$. If for all $n',n''\in S'_n$, $d(n',n'')\leq 2w$, we added a 
station $s'$ with $\N(s')=S'_n\cup \{n\}$ to $S$, if it had not already been added. The cost and 
capacity of $s'$ was found using (\ref{opt12}).

\section{Numerical Experiment}
\label{sec:R}

\subsection{Traffic \& GIS data }

We test our methodology using trip data from the Transportation Tomorrow 
Survey \citep{ttsconduct,ttsdata} covering the 
Greater Golden Horseshoe area of Ontario, Canada. The city of Toronto comprises 16 planning 
districts. We focused on trips made by car in planning district 1 which contains downtown Toronto 
and the surrounding area. 
This dataset contains 90 traffic zones and 78,549 trips over a 24 hour period, broken down into 5 time 
periods: 6:00-9:00, 9:00-15:00, 15:00-19:00, 19:00-24:00, and 24:00-6:00. ESRI shapefiles were 
used to determine node locations as the centroid of each traffic zone. 

\subsection{Estimating $v_t$}
The average trip lengths over each time period are $l=[2.728,2.467,2.661,2.817,2.615]$, and the 
average required time to park or leave a charging station, $p$, was chosen as 10 minutes. We 
calculated $u$ based on information acquired about the bluecar used by Autolib' in Paris, which can 
travel up to 250 km with a recharge time of approximately 4 
hours \cite{auto2010}. The capacity of a pair of parking spaces was then estimated as 
$v=[6366,12874,8512,10570,12794]$ over the 5 time periods, which have been rounded down to the 
nearest even integer.

\subsection{Parking spaces, station cost, and budget}

The cost of installing a pair of parking spaces and the number of pairs of parking spaces were set 
to vary between [1,2] and [1,3], respectively. We took the centroid of all nodes within our dataset 
and considered this point $P$ to be the most expensive and 
dense part of the city area. A node's distance from $P$ determined its cost and parking space 
capacity, with the closest node having a cost of 3 and the furthest having a cost of 1, with 
prices descending linearly with distance. All nodes within $\frac{1}{3}$ of the largest distance 
from the set of nodes to $P$ had a capacity of 1 pair, between $(\frac{1}{3},\frac{2}{3}]$ had a 
capacity of 2 pairs and the remaining nodes had a capacity of 3 pairs. We set our budget 
$b=0.3(c^N)^Tm^N$. The average cost of a pair of parking spaces was $1.67$. This allowed for up to 
39 average priced pairs of parking spaces to be used.

\subsection{Experimental results}
\label{R}
All experiments were done on a Windows 10 Pro 64-bit, Intel Core i7-7820HQ 2.9GHz processor with 8 
GB of RAM computer using Gurobi 8.01. Table 1 presents the results of using BCSA, EE, and 
solving (\ref{opt1}) with $S=N$. An inital attempt at solving EE to optimality was abandoned after 
two days, after which a time limit was set for 12 hours. We observe that BCSA outperformed EE in 
terms of both computation time and solution quality. Using either solution method, we see that 
adding the ability of sharing stations between nodes greatly decreases the number of balanced rides 
not satisfied, with both methods having less than 14\% of the number of unsatisfied rides as $S=N$. 
The set of nodes and stations are plotted in Figure \ref{T2} of Subsection 7.3 of the Appendix, 
where each station was plotted at the centroid of its neighbouring nodes.

\begin{table}[H]\centering
	\def\arraystretch{1.2}
	\begin{tabular}{r c c c c c}
		\hline 
		& BCSA & EE & $S=N$ \\ 
		\hline
		objective value & 6,318 & 6,342  & 46,231	 \\
		computation time (mins)	& 38.63 & 720 & 1.00 \\
		$|S|$ & 85 & 1,764 & 90 \\		
		\hline		
	\end{tabular}
		\label{T0}
	\caption{Comparison of numerical results using BCSA, EE, and $S=N$.}
\end{table}

\section{Conclusion and future research}
\label{s:C}

Vehicle relocation is an operational burden of one-way electric car sharing networks. We have 
developed a novel approach to minimize its requirement by optimizing the location of charging 
stations so as to maximize the balanced flow in the network. By maximizing the balanced flow over all 
charging stations, the trips we expect to satisfy do not depend upon vehicle relocation, as the 
network self-regulates over time. Though we were interested in creating a self-regulating car sharing 
network without the need for vehicle relocation in this paper, there is nothing preventing the use of 
both to further increase the number of trips satisfied. Future research could involve a two-stage 
stochastic optimization model, where charging stations are optimized, and then given a stochastic 
daily network usage, the number of required employees and their deployment is determined. We also see 
our methodology being useful in other one-way sharing programs, such as for bicycles, which suffer 
from the same network imbalance problem.

\section*{Acknowledgements}

This work was supported by Natural Sciences and Engineering Research Council of Canada Discovery 
Grant programs (RGPIN-2015-06164, RGPIN6524-15). The authors would like to thank Mark Ferguson and 
Sean Spears of the McMaster Institute for Transportation \& Logistics for their assistance in 
obtaining the OD trip and GIS location data for the numerical experiment.

\bibliographystyle{plainnat}
\bibliography{EVchargers}

\section{Appendix}

\subsection{Requirement of (5) in (\ref{opt1}) for trip integrity}
\label{App1}

Consider the single period problem of there being two nodes, $n_1$ and $n_2$, and three stations 
$s_1$, $s_2$, and $s_3$. The capacity of each station is 10 (e.g. $v_1=10$, $m_{s_i}=1$ for all $i$, 
and $b=\infty$). Let 

$$OD^1=\begin{bmatrix}
0 & 10 \\
10 & 0 
\end{bmatrix}$$   

and assume $\N(n_1)=\{s_1\}$ and $\N(n_2)=\{s_2,s_3\}$. We see in total there are 20 trips, 10 going 
from $n_1$ to $n_2$, and 10 going to $n_2$ to $n_1$. As $n_2$ has access to two stations, all of its 
incoming and outgoing flow can be accommodated, 5 inflow and 5 outflow going to each, whereas $n_1$ 
has access only to one station, and so can only assign 5 inflow and 5 outflow. The optimal solution 
is 10, with 

$$E^1=\begin{bmatrix}
0 & 5 \\
5 & 0 
\end{bmatrix}$$\\

If we consider now the case where constraint (5) of (\ref{opt1}) is omitted, the optimal solution is 
5, with 

$$E^1=\begin{bmatrix}
5 & 0 \\
0 & 0 
\end{bmatrix}$$\\

Placing 5 in cell (1,1) accounts for the 5 inflow and 5 outflow not assigned by $n_1$, while not 
influencing the flow allocation of $n_2$, allowing it to allocate all of its flow to its 
neighbouring stations, while in reality 5 of $n_1$'s outgoing trips were never assigned to a 
charging station's outward flow. 

\subsection{An example of $AX^t\{\leq,=\} b^t$}

We write $AX^t\{\leq,=\} b^t$ for the example problem described in Subsection \ref{App1}. The curly 
bracketed numbers in the first column indicate the constraint set each row is from.

$$
\begin{matrix}
(1')\\
(1')\\
(1')\\
(2)\\
(2)\\
(3)\\
(3)\\
(4)\\
(4)\\
(4)\\
(5)\\
(5)\\
(5)\\
(5)
\end{matrix}
\begin{bmatrix}
\mathbf{1} & 0 & 0 & 0 & 0 & 0 & 0 & 0 & 0 & 0\\
0 & \mathbf{1} & 0 & 0 & 0 & 0 & 0 & 0 & 0 & 0\\
0 & 0 & \mathbf{1} & 0 & 0 & 0 & 0 & 0 & 0 & 0\\
\mathbf{1} & 0 & 0 & 0 & 0 & 0 & \mathbf{1} & \mathbf{1} & 0 & 0\\
0 & \mathbf{1} & \mathbf{1} & 0 & 0 & 0 & 0 & 0 & \mathbf{1} & \mathbf{1}\\
0 & 0 & 0 & \mathbf{1} & 0 & 0 & \mathbf{1} & 0 & \mathbf{1} & 0\\
0 & 0 & 0 & 0 & \mathbf{1} & \mathbf{1} & 0 & \mathbf{1} & 0 & \mathbf{1}\\
\mathbf{1} & 0 & 0 & \mathbf{-1} & 0 & 0 & 0 & 0 & 0 & 0\\
0 & \mathbf{1} & 0 & 0 & \mathbf{-1} & 0 & 0 & 0 & 0 & 0\\
0 & 0 & \mathbf{1} & 0 & 0 & \mathbf{-1} & 0 & 0 & 0 & 0\\
0 & 0 & 0 & 0 & 0 & 0 & \mathbf{1} & 0 & 0 & 0\\
0 & 0 & 0 & 0 & 0 & 0 & 0 & \mathbf{1} & 0 & 0\\
0 & 0 & 0 & 0 & 0 & 0 & 0 & 0 & \mathbf{1} & 0\\
0 & 0 & 0 & 0 & 0 & 0 & 0 & 0 & 0 & \mathbf{1} 
\end{bmatrix}
\begin{matrix}
\begin{bmatrix}
F^+_{1,1}\\
F^+_{2,2}\\
F^+_{3,2}\\
F^-_{1,1}\\ 
F^-_{2,2}\\
F^-_{3,2}\\
E_{1,1}\\
E_{2,1}\\
E_{1,2}\\
E_{2,2}\\
\end{bmatrix}\\
\\
\\
\\
\\
\end{matrix}
\begin{matrix}
\leq\\
\leq\\
\leq\\
=\\
=\\
=\\
=\\
=\\
=\\
=\\
\leq\\
\leq\\
\leq\\
\leq
\end{matrix}
\begin{bmatrix}
5\\
5\\
5\\
10\\
10\\
10\\
10\\
0\\
0\\
0\\
0\\
10\\
10\\
0
\end{bmatrix}
$$

\newpage

\subsection{Station maps}

\begin{figure}[H]
	\centering
	\begin{tikzpicture}
	\tikzstyle{every node}=[font=\footnotesize]
	\begin{groupplot}[group style={group name=myplot,group size= 1 by 3, 
		vertical sep=0.5cm},height=7.5cm,width=12cm]
	\nextgroupplot[]
	\addplot[only marks, mark options={mark size=2, draw=lightgray, fill=lightgray}]
	table[x=x,y=y]{xy.dat};\label{plots:plot1}
	\addplot[only marks, mark size=1.75, color=black]
	table[x=x,y=y]{stb.dat};\label{plots:plot2}
	\addplot[only marks, mark size=3, color=gray]
	table[x=x,y=y]{cent.dat};\label{plots:plot8}
	\nextgroupplot[]
	\addplot[only marks,  mark options={mark size=2, draw=lightgray, fill=lightgray}]
	table[x=x,y=y]
	{xy.dat};\label{plots:plot1}
	\addplot[only marks, mark size=1.75, color=black]
	table[x=x,y=y]
	{ste.dat};\label{plots:plot2}
	\addplot[only marks, mark size=3, color=gray]
	table[x=x,y=y]
	{cent.dat};\label{plots:plot8}
	\nextgroupplot[]
	\addplot[only marks, mark options={mark size=2, draw=lightgray, fill=lightgray}]
	table[x=x,y=y]
	{xy.dat};\label{plots:plot1}
	\addplot[only marks, mark size=1.75, color=black]
	table[x=x,y=y]
	{stn.dat};\label{plots:plot2}
	\addplot[only marks, mark size=3, color=gray]
	table[x=x,y=y]
	{cent.dat};\label{plots:plot8}
	\end{groupplot}
	\path (myplot c1r1.north west|-current bounding box.north)--
	coordinate(legendpos)
	(myplot c1r1.north east|-current bounding box.north);
	\matrix[matrix of nodes,anchor=south,draw,inner sep=0.2em,draw]at([yshift=1ex]legendpos)
	{	\ref{plots:plot2}& station&[5pt]
		\ref{plots:plot1}& node&[5pt]
		\ref{plots:plot8}& $P$&[5pt]\\};
	\end{tikzpicture}
	\caption{Traffic nodes and station locations for BCSA (top), EE (middle), and $S=N$ (bottom).} 
	\label{T2}
\end{figure}
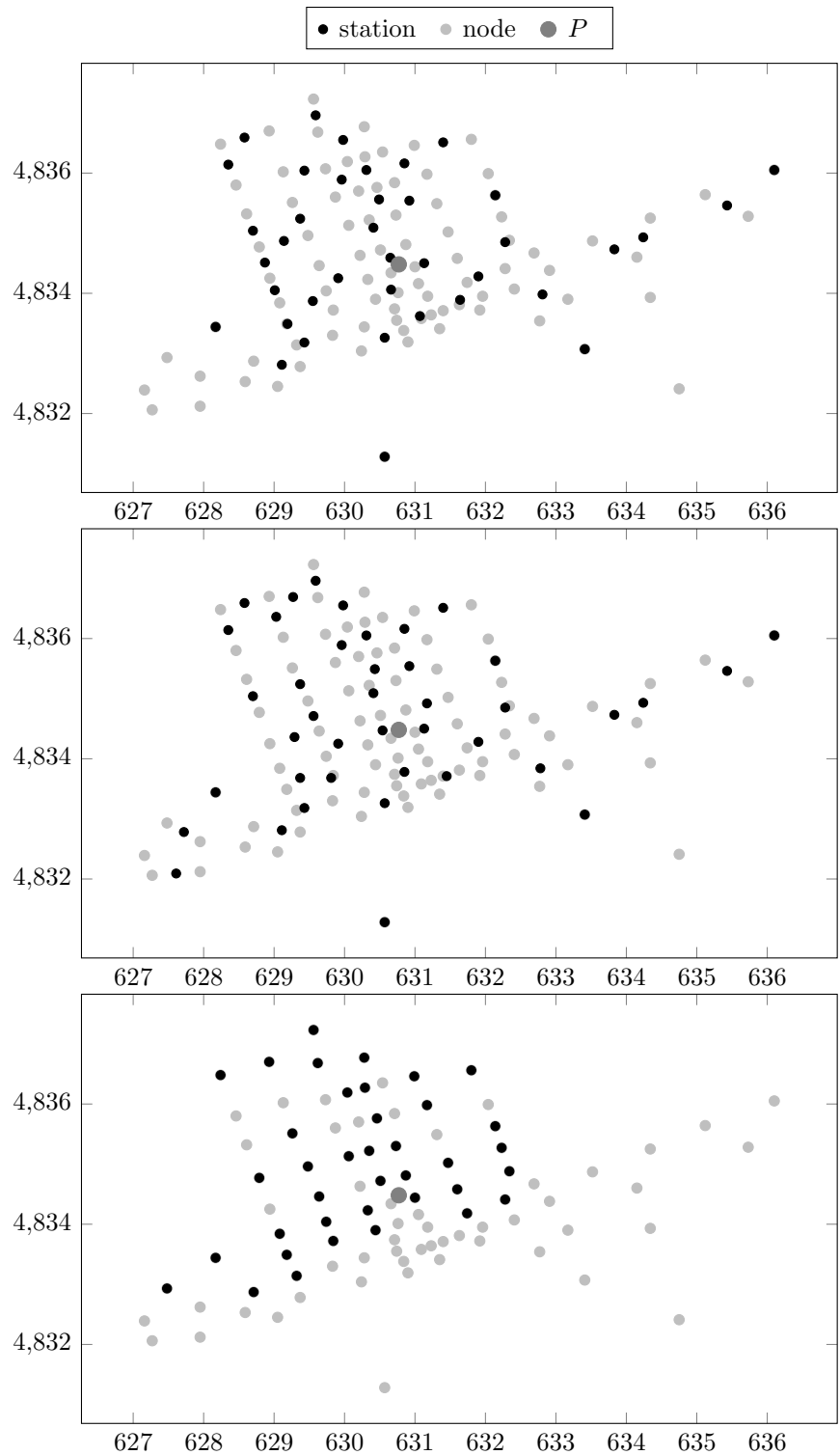

\end{document}